\documentclass[11pt]{amsart} 
\usepackage{amssymb,amsmath,enumerate}

\newtheorem{theorem}{Theorem} 
\newtheorem{corollary}[theorem]{Corollary}
\newtheorem{proposition}[theorem]{Proposition}
\newtheorem{lemma}[theorem]{Lemma}
\newtheorem{conjecture}{Conjecture}
\newtheorem*{example}{Example}
\newtheorem*{remark}{Remark}

\newtheorem*{definition}{Definition}
\newenvironment{defn}{\begin{definition}\rm}{\end{definition}}

\title{Cyclotomic Polytopes and Growth Series of Cyclotomic Lattices}

\author{Matthias Beck and Serkan Ho\c{s}ten}
\address{Department of Mathematics\\
         San Francisco State University\\
         San Francisco, CA 94132\\
         USA}
\email{[beck,serkan]@math.sfsu.edu}

\date{14 February 2006}

\newcommand{\C}{{\mathcal C}}
\renewcommand{\P}{{\mathcal P}}
\renewcommand{\L}{{\mathcal L}}
\newcommand{\R}{{\mathbb R}}
\newcommand{\Z}{{\mathbb Z}}

\newcommand{\0}{{\mathbf 0}}
\newcommand{\1}{{\mathbf 1}}

\newcommand\conv{\operatorname{conv}} 
\newcommand\cone{\operatorname{cone}}

\def\th{^{\text{th}}}
\def\IN{\mathsf{in}}

\newcommand\comment[1]{}
\renewcommand\comment[1]{\textbf{\textsf{[#1]}}} % comment this line out to hide comments

\begin{document}

\begin{abstract} 
The coordination sequence of a lattice $\L$ encodes the word-length function with respect to $M$, a set that generates
$\L$ as a monoid. 
We investigate the coordination sequence of the cyclotomic lattice $\L = \Z[\zeta_m]$, where $\zeta_m$ is a primitive $m\th$ root of unity and where $M$ is the set of all $m\th$ roots of unity. 
We prove several conjectures by Parker %\cite{parkercodeconjectures}  
regarding the structure of the rational generating function of the coordination sequence; this structure depends on the prime factorization of $m$.
Our methods are based on unimodular triangulations of the $m\th$ cyclotomic polytope, the convex hull of the $m$ roots of unity in $\R^{ \phi(m) }$, with respect to a canonically chosen basis of $\L$.
\end{abstract}

\maketitle

%%%%%%%%%%%%%%%%%%%%%%%%%%%%%%%%%%%%%%%%%%%%%%%%%%%%%%%%%%%%%%%%%%%%%%%%%%%

\section{Introduction}

Let $\L \subset \R^d$ be a lattice of rank $r$, and let $M$
be a subset that generates $\L$ as a monoid.  
The \emph{coordination sequence} 
$\left( S (n) \right)_{ n \ge 0 }$ of $(\L, M)$ 
is given by $S(n)$, the number of elements in $\L$ with \emph{word  
length} $n$ with respect to $M$, that is, 
the number of lattice elements that are expressed as a sum from $M$
with a minimal number of $n$ terms \cite{conwaysloanecellstructure}.
The \emph{growth series} $G$ of $(\L, M)$ is the generating function of 
$S (n)$:
\[
  G (x) := \sum_{ n \ge 0 } S (n) \, x^n .
\]
Benson \cite{bensongrowthseries} proved that $G (x) = \frac{ h(x) }{ 
(1-x)^r }$ where $h(x)$, the \emph{coordinator polynomial} of $\L$, 
is a polynomial of degree $\leq  r$. Consequently, 
$S(n)$ is a polynomial of degree $r-1$.
The rationality of $G(x)$ when $\L \cong \Z^r$ is an easy by-product 
of our approach we present below (see also \cite{wagreich}).

Now let $\zeta_m := e^{ 2 \pi i/m }$. We denote by $\Phi_m(x)$ the $m\th$ 
cyclotomic polynomial; its degree is $\phi(m)$, the Euler totient function.
The ring of integers in the cyclotomic field of order $m$, 
$\Z [\zeta_m]$, is a lattice of full rank in $\Z [\zeta_m] 
\otimes_\Z \R \cong \R^{ \phi(m) }$ and hence isomorphic to $\Z^{ \phi(m) }$.
For the remainder of the paper, we let $M$ be the set of all $m\th$ roots of unity, 
and we let $h_m(x)$ be the corresponding coordinator polynomial.  
The study of the coordination sequence of $\Z [\zeta_m]$ with respect to 
$M$ was initiated by Parker, who was motivated by applications to 
error-correcting codes and random walks. 
His article \cite{parkercodeconjectures} includes Kl\o ve's proof 
of the following result, previously conjectured by Parker.

\begin{theorem}[Kl\o ve--Parker] \label{coordinatorprime}
The coordinator polynomial $h_p(x)$ of $\Z[ \zeta_{p}]$, where $p$ is prime, equals
\[
  \Phi_p (x) = x^{ p-1 } + x^{ p-2 } + \dots + 1 \ .
\]
\end{theorem}
Kl\o ve's proof uses a counting argument that relates elements of $\Z[\zeta_p]$ to ordered partitions.

Parker \cite{parkercodeconjectures} offered several conjectures. 
We call a degree-$d$ polynomial $c_d x^d + c_{ d-1 } x^{ d-1 } + \dots + c_0$ \emph{palindromic} if $c_k = c_{ d-k }$.

\begin{conjecture}[Parker] \label{parkerconj1}
The coordinator polynomial $h_m(x)$ of $\Z[\zeta_m]$ equals $\left( h_{ \sqrt m } (x) \right)^{ m/\sqrt m }$, 
where $\sqrt m$ is the squarefree part of $m$.
Furthermore, $h_{ \sqrt m } (x)$ is a palindromic polynomial of degree $\phi(\sqrt m)$.
\end{conjecture}

\begin{conjecture}[Parker] \label{parkerconj2}
The coordinator polynomial of $\Z [\zeta_{2p}]$, where $p$ is an 
odd prime, equals
\[
  h_{ 2p } (x) = 
  \sum_{ k=0 }^{ \frac{ p-3 }{ 2 } } \left( x^k + x^{ p-1-k } \right) \sum_{ j=0 }^{ k } \binom p j + x^{ \frac{ p-1 }{ 2 } } \sum_{ j=0 }^{ \frac{ p-1 }{ 2 } } \binom p j 
  = \sum_{ k=0 }^{ \frac{ p-3 }{ 2 } } \left( x^k + x^{ p-1-k } \right) 
\sum_{ j=0 }^{ k } \binom p j + 2^{p-1} x^{ \frac{ p-1 }{ 2 } }. 
\]
\end{conjecture}

\begin{conjecture}[Parker] \label{parkerconj3}
The coordinator polynomial of $\Z \left[ \zeta_{15} \right]$ equals
\[
  h_{ 15 } (x) = \left( 1 + x^8 \right) + 7 \left( x + x^7 \right) + 28 \left( x^2 + x^6 \right) + 79 \left( x^3 + x^5 \right) + 130 x^4 .
\]
\end{conjecture}

Patras and Sol\'e studied Theorem \ref{coordinatorprime} and 
Parker's conjectures from the viewpoint of Ehrhart polynomials of 
the \emph{cyclotomic polytope} of $\Z[\zeta_m]$ (which we will define 
below). Their article \cite{patrassole} includes an alternate proof of Theorem 
\ref{coordinatorprime} and a computation 
of $h_{ 2p } (x)$ that gave further credence to 
Conjecture \ref{parkerconj2}.

In this paper we prove Conjectures \ref{parkerconj2} and \ref{parkerconj3}, 
and we partly confirm Conjecture \ref{parkerconj1}, in form of the following two theorems.

\begin{theorem}\label{factorizationthm}
For any positive integer $m$, the coordinator polynomial of $\Z[\zeta_m]$ equals $\left( h_{ \sqrt m } (x) \right)^{ m/\sqrt m }$.
\end{theorem}

\begin{theorem}\label{main}
Suppose the positive integer $m$ is one of the following:
\begin{enumerate}[{\rm (i)}]
\item $m = p^\alpha$ where $p$ is prime,
\item $m = p^\alpha q^\beta$ where $p$ and $q$ are distinct primes, or
\item $m= 2^\alpha p^\beta q^\gamma$ where $p$ and $q$ are distinct odd primes.
\end{enumerate}
Then the coordinator polynomial $h_m(x)$ of $\Z[\zeta_m]$ is of the form $h(x)^{ m/\sqrt m }$, where $h(x)$ is the $h$-polynomial of a simplicial polytope, and hence it is palindromic, unimodal, and has nonnegative integer coefficients.
\end{theorem}

Our methods are based on unimodular triangulations of the cyclotomic polytope $\C_m$, 
which we introduce in Section \ref{cycpolsection}. 
We show how one can compute 
$\C_m$ from $\C_{p_1}, \ldots, \C_{p_k}$ where 
$m = p_1^{\alpha_1} \cdots p_k^{\alpha_k}$. 
In Section \ref{hilbseriessection} we study the 
Hilbert series of $\C_m$ and its connection to the growth series 
$G(x)$, and prove Theorem \ref{factorizationthm}. 
We further show that when $m$ is as in one of
the cases of Theorem \ref{main}, the cyclotomic polytope $\C_m$
is totally unimodular. 
In Section \ref{palindromysection} we review toric initial ideals of $\C_m$ and the Dehn--Sommerville relations, and prove Theorem \ref{main}. 
In Section \ref{expcompsection} we compute the face numbers of $\C_p$, $\C_{ 2p }$, 
and $\C_{ 15 }$, and prove Conjectures \ref{parkerconj2} and \ref{parkerconj3}. 
When $p,q,r$ are distinct odd primes then the cyclotomic polytope for $m = pqr$
is {\em not}  totally unimodular. This might be seen as an evidence 
that Conjecture \ref{parkerconj1} may not be true in general. 
In fact, in Section \ref{expcompsection} we present and support the conjecture 
that $h_{ 105 } (x)$ is \emph{not} palindromic.

We would like to point out the recent paper \cite{martinreinercyclotomic} which studies the matroid defined by vertices of the cyclotomic polytope $\C_m$ and its dual matroid, in order to give an upper bound for the number of bases of this matroid. In the cases described in Theorem \ref{main}, this upper bound gives the exact count.
Theorem \ref{duality} below establishes a polytope duality between $\C_m$ and certain multidimensional transportation polytopes and implies that $\C_{pq}$ is simplicial.

%%%%%%%%%%%%%%%%%%%%%%%%%%%%%%%%%%%%%%%%%%%%%%%%%%%%%%%%%%%%%%%%%%%%%%%%%%%

\section{Cyclotomic Polytopes}\label{cycpolsection}

We will now define the \emph{$m\th$ cyclotomic polytope} $\C_m$ associated 
to $\Z[\zeta_m]$. To this end, we will choose a specific lattice basis of $\Z[\zeta_m]$
consisting of certain powers of $\zeta_m$. These powers will correspond to 
the standard unit vectors of $\R^{ \phi(m) }$. The other powers are integer
linear combinations of this basis; hence they are lattice vectors in $\R^{ \phi(m) }$.
The $m\th$ cyclotomic polytope $\C_m$ is the convex hull of all of these 
$m$ lattice points in $\R^{ \phi(m) }$ which correspond to the $m\th$ roots
of unity.
We give this construction first when $m$ is prime, then for a prime power $m$, 
and finally when $m$ is the product of two relatively prime integers. These 
three cases will define $\C_m$ for any positive integer $m$.

When $m = p$ is a prime number
we fix the $\Z$-basis $1, \zeta_p, \zeta_p^2, \ldots, \zeta_p^{p-2}$ of
the lattice $\Z[\zeta_p]$. Since $\zeta_p^{p-1} = 
-\sum_{i=0}^{p-2} \zeta_p^i$,   these $p$ elements form a monoid basis
for $\Z[\zeta_p]$. We identify them with $e_0, e_1, \ldots, e_{p-2}, 
-\sum_{i=0}^{p-2} e_i$ in $\Z^{p-1}$. Hence we obtain:

\begin{proposition}\label{cyc-prime}
The cyclotomic polytope $\C_p \subset \R^{p-1}$, for $p$ prime, is the simplex
\[
\C_p \quad = \quad \conv \left(e_0, \, e_1, \ldots, e_{p-2}, \, 
-\sum_{i=0}^{p-2} e_i \right).
\] 
The only interior lattice point of $\C_p$ is the origin. 
\end{proposition} 

In order to describe $\C_m$ for general $m$ we need two operations on polytopes. 
The first one is the \emph{direct sum} (sometimes called \emph{free sum}; see 
\cite{henkrichtergebertziegler,mcmullenpolconstructions,perlesshephard}). 
Let $P \subset \R^{d_1}$ and $Q \subset \R^{d_2}$ be two polytopes each of which
contains the origin in its interior. Then we define 
\[
  P \circ Q := \conv \left(P \times \0_{d_2}, \0_{d_1} \times Q 
\right) 
\subset \R^{d_1 + d_2}.
\]
Here $\0_d$ denotes the origin in $\R^d$. The polytope
$P \circ Q$ contains $\0_{d_1+d_2}$ in its interior
and its dimension  is the sum of the dimensions of $P$
and $Q$. We denote the $k$-fold direct sum
$P \circ \cdots \circ P$ by $P^{\circ k}$. 

For a prime $p$ and an integer $\alpha \ge 2$, let $\zeta := \zeta_{p^\alpha}$ 
be a primitive $({p^\alpha})\th$ root of unity.
The powers $\zeta^{k+jp^{\alpha-1}}$, where $0 \le k \le p^{\alpha -1}-1$
and $0 \le j \le p-2$, form a lattice
basis of $\Z[\zeta_{p^\alpha}]$, and we will identify them
with the standard unit vectors in $\Z^{\phi(p^\alpha)} = \Z^{ p^{ \alpha -1 } (p-1) } $. 
When we do the identification as  $\zeta^{k+jp^{\alpha-1}} \longleftrightarrow 
e_{k(p-1)+j}$, the cyclotomic polytope $\C_{p^\alpha}$ is 
\[
\C_{p^\alpha} = 
\conv \left(e_{k(p-1)+0}, \, e_{k(p-1)+1}, \, \ldots, \, e_{k(p-1)+p-2}, 
\,
-\sum_{n=0}^{p-2} e_{k(p-1)+n} \,\, : \,\, k = 0, \ldots, p^{\alpha -1}-1
\right).
\] 

\begin{proposition} \label{cyc-primepower}
The cyclotomic polytope $\C_{p^\alpha}$, where $p$ is prime, is equal
to $\C_{p}^{\circ p^{\alpha-1}}$. This polytope is a simplicial
polytope of dimension $\phi(p^\alpha) = p^{\alpha-1}(p-1)$ and the
origin is the only interior lattice point. 
\end{proposition}

\begin{proof} 
As above, let $\zeta := \zeta_{p^\alpha}$ be a primitive
$({p^\alpha})\th$ root of unity. Since $\Phi_p \left( \zeta^{p^{\alpha-1}} \right) = 0$ 
we have 
\[
-1 - \zeta^{p^{\alpha-1}} -\zeta^{2p^{\alpha-1}} - \dots -
\zeta^{(p-2)p^{\alpha-1}} = \zeta^{(p-1)p^{\alpha-1}}. 
\] 
By multiplying this expression with $\zeta^k$ for 
$k=0, \ldots, p^{\alpha-1}-1$ we get
\[
-\zeta^k - \zeta^{k+p^{\alpha-1}} -\zeta^{k+2p^{\alpha-1}} - \dots -
\zeta^{k+(p-2)p^{\alpha-1}} = \zeta^{k+(p-1)p^{\alpha-1}}. 
\] 
The roots of unity that appear on the left-hand side are all
distinct and they are $\zeta^j$ for $j=0, \ldots, (p-1) p^{\alpha-1} - 1$.
This is our chosen lattice basis of $\Z[\zeta]$.

By Proposition \ref{cyc-prime}, $\C_{p^\alpha}$ is precisely $\C_p^{\circ 
p^{\alpha-1}}$. It follows that $\C_{p^\alpha}$ is simplicial 
since $P \circ Q$ is simplicial if $P$ and $Q$ are. 
%Finally, the fact that the origin is the only interior lattice point is clear.
\end{proof}

\begin{example}[The cyclotomic polytope $\C_9$] \rm
To clarify the proof of Proposition \ref{cyc-primepower} we treat the case
$m = 9 = 3^2$. Since $\Phi_3(x) = 1 + x + x^2$ we 
get 
\[ -1 -\zeta^3 = \zeta^6, \quad -\zeta - \zeta^4 = \zeta^7, \quad 
-\zeta^2 - \zeta^5 = \zeta^8 ,
\]
where $\zeta = \zeta_9$ is a primitive $9\th$ root of unity. So 
\[
\C_9 = \conv \left(e_0, \, e_1, \,- e_0-e_1, \,\,\, e_2, 
\, e_3, \, -e_2-e_3,  \,\,\, e_4, \, e_5, \, -e_4-e_5\right) ,
\] 
and this is exactly $\C_3 \circ \C_3 \circ \C_3$.
\end{example}

We now recursively construct a lattice basis for $\Z[\zeta_m]$
and the cyclotomic polytope $\C_m$, by factoring $m = m_1 m_2$,
where $m_1$ and $m_2$ are relatively prime, and assuming that 
bases for $\Z[\zeta_{m_1}]$ and $\Z[\zeta_{m_2}]$ and the cyclotomic
polytopes $\C_{ m_1 }$ and $\C_{ m_2 }$ in these bases are already 
constructed.

So assume that $\omega_1, \ldots, \omega_{\phi(m_1)}$ form a $\Z$-basis 
of $\Z[\zeta_{m_1}]$, and together with $\omega_{\phi(m_1)+1},  \ldots,
\omega_{m_1}$ they form a monoid basis. Then
\[
  \C_{m_1} = \conv \left( e_1, \ldots, e_{\phi(m_1)}, v_{\phi(m_1)+1}, \ldots, v_{m_1} \right) \subset \R^{\phi(m_1)} .
\]
Similarly, we assume
that $\rho_1, \ldots, \rho_{\phi(m_2)}$ form a $\Z$-basis 
of $\Z[\zeta_{m_2}]$, and together with $\rho_{\phi(m_2)+1}, \ldots, 
\rho_{m_2}$ they form a monoid basis. Now
\[
  \C_{m_2} = \conv \left( f_1, \ldots, f_{\phi(m_2)}, w_{\phi(m_2)+1}, \ldots, w_{m_2} \right) \subset \R^{\phi(m_2)} .
\]
For the cyclotomic lattice
$\Z[\zeta_m]$ the set of $m\th$ roots $\{\omega_i \rho_j \,\, : \,\,
1 \leq i \leq \phi(m_1), \,\, 1 \leq j \leq \phi(m_2) \}$
is a basis, and the pairwise product of the lattice points in 
$\R^{ \phi(m) }$ corresponding to \emph{all} of the $m_1\th$
and $m_2\th$ roots is a monoid basis of $\Z[\zeta_m]$. 
We define the cyclotomic polytope $\C_m$ to be the convex hull 
of the vectors in this monoid basis.

This construction motivates the second polytope operation we need,
namely the tensor product of two polytopes.
Let $P \subset \R^{d_1}$ and $Q \subset \R^{d_2}$ be two polytopes
with vertices $v_1, \ldots, v_s$ and $w_1, \ldots, w_t$, respectively.
Then $P \otimes Q \subset \R^{d_1d_2}$ is the polytope
\[
P \otimes Q := \conv \left(v_i \otimes w_j \, : \, 
1 \leq i \leq s, \,\,\, 1 \leq j \leq t\right).
\]
Our construction of $\C_m$ immediately implies:

\begin{proposition} \label{cyc-product}
Let $m = m_1m_2$ where $m_1, m_2 > 1$ are relatively prime. 
Then the cyclotomic polytope $\C_m$ is equal to $\C_{m_1} \otimes 
\C_{m_2}$. 
\end{proposition}

Propostions \ref{cyc-prime}, \ref{cyc-primepower}, and \ref{cyc-product}
allow us to describe the cyclotomic polytope $\C_m$ for any positive 
integer $m$. Moreover, $\C_m$ is determined by 
$\C_{ \sqrt m }$.

\begin{theorem} \label{cyc-thm}
The cyclotomic polytope $\C_m$ is equal
to $\C_{\sqrt{m}}^{\circ (m/\sqrt{m})}$ where $\sqrt{m}$ is 
the squarefree part of $m$.
\end{theorem}

\begin{proof} Let $m = p_1^{\alpha_1} \cdots p_n^{\alpha_n}$. 
Propositions \ref{cyc-primepower} and \ref{cyc-product}
imply that $\C_m = \bigotimes_{i=1}^n \C_{p_i}^{\circ p^{\alpha_i-1}}$. 
If we let $V_i$ be the matrix whose columns are the vertices of 
$\C_{p_i}$, then the polytope  $\C_{p_i}^{\circ p^{\alpha_i-1}}$ is
equal to $\conv \left( V_i \otimes I_{p^{\alpha_i-1}} \right)$, 
where $I_k$ denotes a $k \times k$ identity matrix.
Therefore 
\[
  \C_m 
  = \conv \left( \bigotimes_{i=1}^n (V_i \otimes I_{p^{\alpha_i-1}}) \right) 
  = \conv \left( \bigotimes_{i=1}^n V_i \otimes \bigotimes_{i=1}^n I_{p^{\alpha_i-1}} \right) .
\]
The last expression is precisely $\C_{\sqrt{m}}^{\circ (m/\sqrt{m})}$.
\end{proof}

\begin{lemma} \label{zeroplusminusone}
Let $m = p_1 p_2 \cdots p_k$,
where $p_1, p_2, \dots, p_k$ are distinct primes.
Then the vertices of $\C_m$ have coordinates in $\{0, +1, -1\}$,
and they are precisely the tensor products of the vertices
of $\C_{p_1}, \ldots, \C_{p_k}$. The only other lattice point in
$\C_m$ is the origin, which is in the interior of $\C_m$.
\end{lemma}

\begin{proof}
The vertices of $\C_m$ have coordinates in $\{0, +1, -1\}$ by construction.
It is clear that none of the tensor products of the vertices
of $\C_{p_1}, \ldots, \C_{p_k}$ is a convex combination of the others, so 
they are all vertices.
It remains to prove that the origin is the only other lattice point in $\C_m$.
%Theorem \ref{cyc-thm} reduces our discussion to the case where $m$ is squarefree.
If $m$ is prime, Proposition \ref{cyc-prime} says that there is no other lattice point aside from the origin and the vertices.
Now let $m = n p$ where $n$ is squarefree, not divisible by $p$, and $p$ is prime.
By induction, $\C_n$ does not contain any lattice point other than the origin and its vertices.
Let $A_m$ be the matrix whose columns are the vertices of $\C_m$, then
\[
A_m = \left[ \begin{array}{ccccc}
A_n &     &  &  &  -A_n \\
    & A_n &  &  &  -A_n \\
    &     &\ddots&  & \vdots \\
    &     &   & A_n& -A_n 
\end{array} \right].
\]
(This matrix has $p$ column blocks and $p-1$ row blocks.)
Now suppose there is a nonzero lattice point $u \in \C_m$ that is a nontrivial 
convex combination of the vertices. The point $u$ is a $0, \pm 1$ vector, and we may assume that it has first coordinate $1$. Then $u$ has to be a convex combination of vertices of $\C_m$ that have $1$ as the first coordinate.
This means that $u$ is a convex combination of such vectors coming from the 
first $A_n$-block of $A_m$ and the top $(-A_n)$-block.
By looking at the coordinates of the first row of the second block of $A_n$'s 
in $A_m$, we see that the corresponding coordinate of $u$ is strictly between 
$0$ and $1$ if in the convex combination vectors from both the first $A_n$ 
and the top $(-A_n)$-block were used. 
%Hence $u$ came from a convex combination of just $A_n$, 
Hence $u$ must be a convex combination of the first $n$ columns of $A_m$,
which contradicts our induction hypothesis. 
\end{proof}

\begin{corollary} \label{onepoint}
The cyclotomic polytope $\C_m$ is a 
$\{0,+1,-1\}$-polytope with only one lattice point other than
its vertices. This lattice point is the origin and it is in the
interior of $\C_m$.
\end{corollary}
\begin{proof} This follows from Theorem \ref{cyc-thm} and Lemma \ref{zeroplusminusone}.
\end{proof}

We will return to the combinatorial structure of $\C_m$ in 
Section \ref{palindromysection}.

%%%%%%%%%%%%%%%%%%%%%%%%%%%%%%%%%%%%%%%%%%%%%%%%%%%%%%%%%%%%%%%%%%%%%%%%%%%

\section{Hilbert Series and Unimodular Triangulations}\label{hilbseriessection}

Let $\L \subset \Z^d$ be a lattice, let $M$ be a minimal set of monoid
generators, and let $K$ be an arbitrary field. 
The monoid (semigroup) algebra $K[M']$, where $M' = \{(u,1)\, : 
\, u \in M \cup \{0\}\}$, is a finitely generated graded $K$-algebra where 
each monomial in $K[M']$ corresponds to $(v,k)$ where 
$v = \sum_{u_i \in M \cup \{0\}} n_iu_i$ with nonnegative integer coefficients $n_i$ such that $\sum n_i = 
k$. Such an element has degree $k$ in $K[M']$. In this setting
the {\em Hilbert series} of $K[M']$ is
\[
H(K[M']; x) \, \, := \, \, \sum_{k \geq 0} \dim_K \left( K[M']_k \right) x^k 
,
\]
where $K[M']_k$ denotes the vector space of elements of degree $k$ in this graded algebra.
It is a standard result of commutative algebra that
\[
H(K[M']; x) \, = \, \frac{h(x)}{(1-x)^{d+1}} 
\] 
where $h(x)$ is a polynomial of degree at most $d$ 
\cite{atiyahmacdonald, brunsherzog}.
When $\L \cong \Z^d$, it is clear that the number of elements in $\L$ of 
length 
exactly $k$ (with respect to $M$) is equal to 
$\dim_K(K[M']_k) - \dim_K(K[M']_{k-1})$, and therefore 
the growth series is 
\[
G(x) \, = \, (1-x)H(K[M']; x) \, = \, \frac{h(x)}{(1-x)^d} \ ,
\]
and this reproves the rationality of $G(x)$ in this case.

\begin{lemma}\label{monoidmult}
Let $K[M']$ be the monoid algebra corresponding to the cyclotomic polytope $\C_m$.
Let $N$ be the lattice points in $\C_m \circ \C_m$ and $K[N']$ be the corresponding
monoid algebra. Then
\[
  H(K[N']; x) = H(K[M']; x) \cdot H(K[M']; x) \ .
\]
\end{lemma}

\begin{proof}
This follows from
\[
  \dim_K \left( K[N']_k \right) = \sum_{ s+t = k } \dim_K \left( 
K[M']_s \right) \dim_K \left( K[M']_t \right) .
\]
\end{proof}

The statement of Theorem \ref{factorizationthm}, namely that
$h_m (x) = \left( h_{ \sqrt m } (x) \right)^{ m/\sqrt m }$,
follows now immediately from Theorem \ref{cyc-thm} and Lemma \ref{monoidmult}:

\begin{proof}[Proof of Theorem \ref{factorizationthm}]
Let $N = \C_m \cap \Z^{ \phi(m) } $ 
and $M = \C_{ \sqrt m } \cap \Z^{ \phi(\sqrt m) } $, and let
$K[N']$ and $K[M']$ be the corresponding monoid algebras.
Theorem \ref{cyc-thm} implies $\C_m = \C_{\sqrt{m}}^{\circ (m/\sqrt{m})}$
and Lemma \ref{monoidmult} implies that
\[
  H(K[N']; x) = \left( H(K[M']; x) \right)^{ m/\sqrt m } .
\]
This means that $h_m (x) = \left( h_{ \sqrt m } (x) \right)^{ m/\sqrt m }$.
\end{proof}

Now let $\P_M \subset \R^d$ be the convex hull of $M$. 
Suppose that the set of lattice points $\P_M \cap \Z^d$ is equal
to $M \cup \{0\}$. Corollary \ref{onepoint} implies that $\C_m$ has this
property. In general, the monoid generated by $M'$ and the monoid
of the lattice points in the cone generated by $M'$ are not equal. In the
case of the equality we call $M'$, $\P_M$, and $K[M']$ {\em normal}. 
We give a necessary condition for the normality of these objects below.
Note that when $\P_M$ is normal then the set of lattice points in
$\cone(M') \cap \{x: \, x_{d+1} = k\}$ is in bijection with the set
of lattice points in $k \P_M$, the $k\th$ dilate of $\P_M$.

A simplex with vertices $v_0, v_1, \dots, v_d \subset \Z^d$ is
\emph{unimodular} if $\left\{ v_1 - v_0, v_2 - v_0, \dots, v_d - v_0
\right\}$ generates $\Z^d$. This is equivalent to $|\det(v_1 - v_0, v_2 -
v_0, \dots, v_d - v_0)| = 1$.  A \emph{unimodular triangulation} of a
polytope $\P$ is a triangulation into unimodular simplices with vertices
in $\P \cap \Z^d$. 

\begin{lemma}\label{unimodular-implies-normal} If
$\P_M$ has a unimodular triangulation then $\P_M$ is normal. 
\end{lemma}

\begin{proof} For each unimodular simplex $\sigma = \{v_0, v_1, \ldots,
v_r\}$ in this triangulation we consider 
\[
\cone(\sigma) = \cone \left(
\left(\begin{array}{c} 1 \\ v_0 \end{array} \right), \,\,
\left(\begin{array}{c} 1 \\ v_1 \end{array} \right), \,\, \ldots, \,\,
\left(\begin{array}{c} 1 \\ v_r \end{array} \right) \right) 
\] 
These cones cover $\cone(M')$ and the absolute value of the 
determinant of their generators is one. If $z \in \cone(M')$ then
$z$ is in one of the $\cone(\sigma)$, and if it has integer coordinates, 
Cramer's Rule implies that $z$ is a nonnegative integer linear 
combination of the generators of this cone. This shows that $M'$ is 
normal.
\end{proof} 

\begin{defn} 
A matrix with $0, \, +1, \, -1$ entries is called {\em totally 
unimodular} if every square submatrix has determinant $0, \, +1$, or $-1$. 
We say that a polytope $\P$ is {\em totally unimodular} if the matrix whose columns are
the lattice points in $\P$ is totally unimodular.
\end{defn}

If $\P$ contains the origin in its interior and is totally unimodular, then any 
triangulation of $\P$ that is a cone with apex the origin is unimodular. We will show that the polytope
$\C_m$, where $m$ is an integer described in Theorem \ref{main}, is 
totally unimodular.  We will use the following characterization
of totally unimodular matrices.

\begin{theorem} 
\label{thm-columns}
\cite[Theorem 19.3]{schrijver}
A matrix $A$ with $0,\, +1, \, -1$ entries is totally unimodular
if and only if each collection of columns of $A$ can be split into
two parts so that the sum of the columns in one part minus the
sum of the columns in the other part is a vector with entries
only $0$, $+1$, and $-1$.  
\end{theorem}

\begin{theorem} \label{unimodular}
The polytope $\C_m$ is totally unimodular for all $m$ described in 
Theorem \ref{main}.
\end{theorem}

\begin{proof}
Given $\C_m$, we let $A_m$ be the matrix of its vertices.
By Corollary \ref{onepoint} this matrix has $0, \, +1, \, -1$
entries. 
The matrix 
$\left[ \begin{array}{cc} B & 0 \\ 0 & C \end{array} \right]$ 
is totally unimodular if and only if $B$ and $C$ are.
By Theorem \ref{cyc-thm}, we may assume that $m$ is squarefree,
that is, $m$ is prime, the product of two primes, or $m=2pq$, where $p$ and $q$ are odd primes.

In the first case  $A_p \, = \, \left[I_{p-1} \,\, -\!\1\right]$ where
$\1$ is a column of all ones. This matrix is clearly totally unimodular.
% since given any subset of the columns we put them only in the first part as referred to in Theorem \ref{thm-columns}. 

The matrix for the case $m=pq$ is  
\[
A_{pq} = \left[ \begin{array}{ccccc}
A_p &     &  &  &  -A_p \\
    & A_p &  &  &  -A_p \\
    &     &\ddots&  & \vdots \\
    &     &   & A_p& -A_p 
\end{array} \right].
\]
Now we use Theorem \ref{thm-columns} and split the columns of $A_{ pq }$
into two parts. Given a subset of the columns
of $A_{pq}$ we put all the columns in the last block into the first
part. The sum of these columns is a vector with entries either 
$0$ and $-1$ only, or $0$ and $+1$ only, depending on 
whether the last column of this block (a $+ \1$) is included 
or not. We treat the second case, and the first case can be
dealt with similarly. We put all the columns that involve 
$-\1$ also in the first part. Now the sum of all these columns
is a vector with $0, \, +1$, and $-1$ only. The remaining
columns are columns of $I_{(p-1)(q-1)}$, and we can arrange them
to be put in the two parts so that the resulting vector
has only $0$, $+1$, and $-1$ entries.

Finally, the matrix $A_{2pq}$ equals $\left[A_{pq} \ -\!A_{pq}\right]$, 
and we immediately conclude that $A_{2pq}$ is also totally unimodular.
\end{proof}

\begin{remark} \rm
Total unimodularity breaks down already in the case of $\C_{ 3pq }$, where $p$ and $q$ are distinct primes $>3$. Here
\[
A_{3pq} = \left[ \begin{array}{ccc}
A_{pq} &        & -A_{pq} \\
       & A_{pq} & -A_{pq}
\end{array} \right] ,
\]
and the columns
\[
  \left( \begin{array}{c}
  \0_{ (p-1)(q-1) } \\
  \1_{ (p-1)(q-1) } 
  \end{array} \right) , \quad
  \left( \begin{array}{c}
  1 \\ \0_{ p-2 } \\ \vdots \\ 1 \\ \0_{ p-2 }
  \end{array} \right) ,
  \quad \text{ and } \quad
  \left( \begin{array}{c}
  \1_{ p-1 } \\ \0_{ (p-1)(q-2) } \\ \1_{ p-1 } \\ \0_{ (p-1)(q-2) }
  \end{array} \right)
\]
violate the condition of Theorem \ref{thm-columns}.
When $m=pqr$ for primes $p,q,r >3$, the polytope $\C_m$ is also
not totally unimodular. This follows from the non-normality of
the monoid algebra of the three-dimensional $(p-1) \times (q-1) 
\times (r-1)$ transportation polytope \cite[p.~77]{vlach}.
Hence $\C_m$ is not totally unimodular when $m$ is
divisible by three or more odd primes.   
\end{remark}
%%%%%%%%%%%%%%%%%%%%%%%%%%%%%%%%%%%%%%%%%%%%%%%%%%%%%%%%%%%%%%%%%%%%%%%%%%%

\section{Palindromy}\label{palindromysection}

The monoid algebra $K[M']$ is a finitely generated graded $K$-algebra, and
hence is isomorphic to $K \left[ x_1, \dots, x_n \right] / I_M$
where $n=|M'|$ and $I_M$ is a homogeneous {\em toric ideal} \cite{GB+CP}. 
For the results in this section we need the notion of 
{\em initial ideals}.

In the polynomial ring $R = K\left[ x_1, \dots, x_n \right]$, 
we abbreviate the monomial $x_1^{u_1} \cdots x_n^{u_n}$ by $x^u$.
A \emph{term order}
$\prec$ is a well ordering of all the monomials in $R$
(with the minimum element $x^0 = 1$) that is
compatible with multiplication; that is, $x^u
\prec x^v$ implies that $x^w x^u \prec
x^w x^v$ for any monomial $x^w$.
Given a nonzero polynomial $f$ and a term order $\prec$,
we let $\IN_{\prec}(f)$, the {\em initial term} of $f$,
be the largest monomial of $f$
with respect to $\prec$.  If $I$ is an ideal,
the {\em initial ideal} of $I$ with respect
to the term order $\prec$ is the monomial ideal
generated by all the initial terms of polynomials in $I$:
                                                                                    
\[
\IN_{\prec}(I) = \langle \IN_{\prec}(f) \,\, : \,\, f \in I \rangle.
\]                                                                                    
                                                                                    
\begin{proposition} 
\label{Hilb-equality}
\cite[Chapter 15]{eisenbud} 
Let $I$ be a homogeneous ideal in $R$ 
and $\prec$ any 
term
order. Then for all $k$, $\left( R/I \right)_k$ and $\left( R/\IN_\prec(I) 
\right)_k$ are isomorphic $K$-vector spaces, and therefore 
$H(R/I; x) = H(R/\IN_\prec(I); x)$.
\end{proposition}
The following result follows from Corollary 8.4 and Corollary 8.9 in \cite{GB+CP}. The regular triangulation $\Delta_\prec$ of $M$ is obtained by 
lifting each point in $M$ by $\omega_i$ where $\omega = 
(\omega_1, \ldots, \omega_n)$ is a weight vector so
that $\IN_\prec(I_M) = \IN_\omega(I_M)$, and then by taking the convex hull 
of these lifted points. The facets of the lower hull of this convex
hull form $\Delta_\prec$; see \cite[Chapter 8]{GB+CP} for more details. 
\begin{theorem} \label{reg-triang}
Let $K[M']$ be a monoid algebra and $I_M$ be the corresponding
toric ideal. The initial ideal $\IN_\prec(I_M)$ is squarefree 
if and only if  the regular triangulation $\Delta_\prec$ of $M$ induced
by $\prec$ is unimodular. In this case $\IN_\prec(I_M)$ is the
Stanley-Reisner ideal of $\Delta_\prec$ viewed as a simplicial
complex.
\end{theorem}

The \emph{$f$-vector} $\left( f_{ -1 } , f_0, \dots, f_{d-1} \right)$ of the $d$-polytope $\P$ consists of its face numbers, so $f_{ -1 } = 1$ (corresponding to the empty face), $f_0$ is the number of vertices of $\P$, $f_1$ the number of edges, and so on, up to $f_{ d-1 }$, the number of facets. Closely related is the \emph{$h$-vector} $ \left( h_0, h_1, \dots, h_d \right)$ of $\P$, defined through
\[
  \sum_{ j=0 }^d h_j \, x^j = \sum_{ k=0 }^{ d } f_{ k-1 } \, (x-1)^{ d-k } \ .
\]
The left-hand side is the \emph{$h$-polynomial} of $\P$.
Explicitly, $h_j$ is given by
\[
  h_j = \sum_{ k=0 }^j (-1)^{ j-k } \binom{ d-k }{ j-k } f_{ k-1 } \ .
\]
The famous \emph{Dehn--Sommerville Relations} assert that, for a simplicial polytope, $h_j = h_{ d-j }$.

We are now ready to prove Theorem \ref{main}, namely that the coordinator polynomial of $\Z[\zeta_m]$, when $m$ is divisible by at most two odd primes, is of the form $h(x)^{ m/\sqrt m }$, where $h(x)$ is the $h$-polynomial of a simplicial polytope. This implies that $h(x)$ is palindromic, unimodal, and has nonnegative integer coefficients.

\begin{proof}[Proof of Theorem \ref{main}] 
Theorem \ref{factorizationthm} reduces the discussion to the case
when $m$ is squarefree.
By Theorem \ref{unimodular} the polytope $\C_m$ is totally 
unimodular. Corollary 
\ref{onepoint} implies that any triangulation of $\C_m$
induced by a triangulation of its boundary (by coning
over the boundary triangulation using the origin as the apex) 
is unimodular. Now we can use a {\em pulling} ({\em reverse lexicographic} 
\cite[Chapter 8]{GB+CP})
triangulation of the boundary of $\C_m$ to obtain such 
a unimodular {\em regular} triangulation $\Delta_\prec$. We note that
this boundary triangulation is the boundary of a simplicial polytope $Q_\prec$
of the same dimension as $\C_m$ (see \cite{stanleypaper} 
for the construction of such a triangulation). 
By Theorem \ref{reg-triang} 
the initial ideal $\IN_\prec(I_M)$ is squarefree and it is the Stanley-Reisner
ideal of $\Delta_\prec$. Proposition  \ref{Hilb-equality} 
implies that $H(K[M']; x) = H(R/\IN_\prec(I_M)); x)$ and
this rational function's numerator is the $h$-polynomial
of $Q_\prec$ \cite[Theorem II.1.4]{stanleybook}. Now the Dehn--Sommerville 
relations 
imply palindromy of the numerator of the growth series. 
Unimodality and nonnegativity follow from \cite[Theorem III.1.1]{stanleybook}.
\end{proof}

\begin{corollary}\label{maincor}
If $\C_m$ is a simplicial polytope then the coordinator polynomial of 
$\Z[\zeta_m]$ equals the $h$-polynomial of~$\C_m$.
\end{corollary}

\begin{proof}
As in the above proof, we use the unimodular regular triangulation $\Delta_\prec$ obtained by coning over the boundary of $\C_m$ using the origin as the apex.
\end{proof}

\begin{example} \rm (The coordinator polynomial of $\Z \left[ \zeta_{180} 
\right] $). 
Since $m = 180 = 2^2 \cdot 3^2 \cdot 5$ the coordinator 
polynomial is of the form $p(x)^6$ where $p(x)$ is the
coordinator polynomial of $\Z[\zeta_{30}]$. Using
the software {\tt 4ti2} \cite{4ti2} we computed 
a reverse lexicographic initial ideal of the toric 
ideal $I_M$ where $M = \C_{30} \cap \Z^8$. This initial
ideal has $615$ squarefree minimal generators. Then we used
the computer algebra system {\tt CoCoA} \cite{cocoa}
to compute the Hilbert series from this monomial 
ideal to obtain
\[
p(x) \, = \, x^8 + 22x^7 + 208x^6 + 874x^5 + 1480x^4 + 874x^3 + 208x^2 + 22x + 
1. 
\]
\end{example}

\begin{remark} \rm
The polytope $\C_m$ is not simplicial in general. 
For example, when $m=30$ the polytope $\C_{ 30 }$ is a
non-simplicial polytope of dimension $8$ with $810$ facets.
This polytope has two types of facets: $450$ of them are 
simplicial, and the rest of them are facets with $10$ vertices.
Proposition \ref{cyc-prime}, Proposition \ref{cyc-primepower}, and Theorem
\ref{cyc-thm} together with Proposition \ref{chapman} below
imply that the other candidates for non-simplicial $\C_m$ for $m<30$ are
$m=15$ and $m=21$. However, in these cases the two polytopes are simplical;
$\C_{15}$ has $360$ facets and $\C_{21}$ has $4410$ facets. Hence
$\C_{30}$ is the smallest non-simplicial cyclotomic polytope. We have also
checked that $\C_{33}$ and $\C_{35}$ are simplicial
with $554400$ and $1134000$ facets, respectively. This led us to the
following result whose proof was suggested by Robin Chapman 
\cite{robinchapman}.
\end{remark}

\begin{proposition} \label{chapman}
The cyclotomic polytope $\C_{ pq }$, where $p$ and $q$ are prime, is 
simplicial.
\end{proposition}

The result follows from a polytope duality between $\C_m$ and certain
multidimensional transportation polytopes. We first introduce these
polytopes. Let $p_1, \ldots, p_k$ be positive integers. A 
\emph{multidimensional table} is a $p_1 \times \cdots \times p_k$ array
of real numbers. We will denote the entries of such a table
${\bf x}$ by $x_{i_1 \, \ldots \, i_k}$. Now suppose 
for each $i=1,\ldots,k$ there is a $(k-1)$-dimensional table ${\bf b^i}$
of size $p_1 \times \cdots \times \widehat{p_i}
\times \cdots \times p_k$ with nonnegative real entries.
Then 
\[
P \left( {\bf b^1}, \ldots, {\bf b^k} \right) \, := \, \left\{ {\bf x} \in \R_{ \ge 0 }^{\prod p_i} 
\, : \, 
\begin{array}{ll}
\sum_{j=1}^{p_1} x_{j \, i_2 \ldots i_k} = {b^1}_{i_2 \ldots i_k} 
&\forall \,\, 
(i_2, \ldots, i_k) \, \in \, [p_2] \times \cdots \times [p_k]\\
\sum_{j=1}^{p_2} x_{i_1 \, j \, i_3 \ldots i_k} = {b^2}_{i_1 i_3 \ldots 
i_k} 
&\forall \,\, (i_1, i_3, \ldots, i_k) \, \in \, [p_1] \times 
[p_3] \cdots \times [p_k]\\
&\vdots \\
\sum_{j=1}^{p_k} x_{i_1 \ldots i_{k-1}j} = {b^k}_{i_1 \ldots 
i_{k-1}}
&\forall \,\, (i_1, \ldots, i_{k-1}) \, \in \, [p_1] \times  
\cdots \times [p_{k-1}]
\end{array}
\right\}
\] 
is a \emph{multidimensional transportation polytope} defined by
the tables ${\bf b^1}, \ldots, {\bf b^k}$. We will be concerned with a 
very particular type of transportation polytopes, namely, given
integers $p_1, \ldots, p_k$ we let ${\bf b^i}$ be the table
all whose entries are equal to $p_i$. Such a transportation
polytope will be denoted by $P(p_1, \ldots, p_k)$. 
For instance, when $k=2$, $P(p_1, p_2)$ is the ``usual" transportation 
polytope consisting of nonnegative $p_1 \times p_2$ matrices with all
row sums equal to $p_2$ and all column sums equal to $p_1$.
Now we can state the duality theorem.

\begin{theorem} \label{duality}
Let $m = p_1p_2\cdots p_k$,
where $p_1, p_2, \dots, p_k$ are distinct primes. Then the
cyclotomic polytope  $\C_m$ and the transportation polytope
$P(p_1, \ldots,p_k)$ are dual to each other. 
\end{theorem}  

\begin{proof}
We will show that the face lattice of $\C_m$ and $P(p_1, \ldots, p_k)$
are dual to each other. First we show that there is a bijection between
the facets of $\C_m$ and the vertices of $P(p_1, \ldots, p_k)$.
Each facet of $\C_m$ is defined by a linear form 
$f(x) = 1$. Now let ${\bf y}$ be the  $p_1 \times \cdots \times p_k$ 
table where 
$y_{i_1 \ldots i_k} = f(v_{i_1}^1 \otimes v_{i_2}^2  \otimes \cdots 
\otimes v_{i_k}^k)$ where $v_{i_j}^j$ is a vertex of $\C_{p_j}$. 
The entries of ${\bf y}$ are at most $1$, and those entries
that are equal to $1$ are in bijection with the vertices on
the facet defined by $f(x)=1$. Since the sum of the vertices
of each $C_{p_j}$ is the origin, we conclude that 
\begin{align*}
\sum_{j=1}^{p_1} y_{j \, i_2 \ldots i_k} &= 0 
\qquad \forall \,\, 
(i_2, \ldots, i_k) \, \in \, [p_2] \times \cdots \times [p_k] \\
\sum_{j=1}^{p_2} y_{i_1 \, j \, i_3 \ldots i_k} &= 0
\qquad \forall \,\, (i_1, i_3, \ldots, i_k) \, \in \, [p_1] \times  [p_3] \cdots \times [p_k] \\
&\vdots \\
\sum_{j=1}^{p_k} y_{i_1 \ldots i_{k-1}j} &= 0 
\qquad \forall \,\, (i_1, \ldots, i_{k-1}) \, \in \, [p_1] \times  
\cdots \times [p_{k-1}].
\end{align*}
Now we define a new table ${\bf x}$ where 
$x_{i_1 \cdots i_k} = 1 - y_{i_1 \cdots i_k}$. This table is a 
point of the transportation polytope $P(p_1, \ldots, p_k)$. 
A facet $F$ of $\C_m$ corresponds to a table in this 
transportation polytope whose zero entries are in bijection 
with the vertices of $\C_m$ that are on $F$. On the other 
hand, a vertex of $P(p_1, \ldots, p_k)$ is defined by setting
some of the entries to zero. This implies that ${\bf x}$ has to be
a vertex of $P(p_1, \ldots, p_k)$, since otherwise
there would be a vertex with more zero entries which in turn
give more vertices of $\C_m$ incident to $F$. This contradiction
shows the bijection between the facets of $\C_m$ and the vertices
of $P(p_1, \ldots, p_k)$. 

To extend this bijection to all 
faces we make the following observation: If $F$ is a face of 
$\C_m$ that is the intersection of the facets
$F_1, \ldots, F_t$ defined by the linear forms $f_1, \ldots, f_t$,
then $f = \sum_{i=1}^t \lambda_i f_i$ such that  $\lambda_i > 0$
and $\sum_{i=1}^t \lambda_i = 1$ can be taken
to be a supporting hyperplane of $F$ which is not a 
supporting hyperplane of faces strictly containing $F$. 
Such an $f$ gives rise to 
${\bf x} = \sum_{i=1}^t \lambda_i {\bf x^i}$ where the ${\bf x^i}$'s
are the vertices of $P(p_1, \ldots, p_k)$
corresponding to the facets $F_1, \ldots, F_t$. 
All such ${\bf x}$ form the relative interior of a face $G$ of 
$P(p_1, \ldots, p_k)$ defined by setting 
those entries  of ${\bf x}$ corresponding to the vertices on
$F$ equal to zero. The vertices of $G$ 
are precisely ${\bf x^1}, \ldots, {\bf x^t}$, since any extra vertex
will translate into one more vertex on $F$.  
\end{proof}

\begin{proof}[Proof of Proposition \ref{chapman}]
The points in $P(p,q)$ are in bijection with nonnegative edge
assignments of the complete bipartite 
graph $K_{p,q}$ with the node partition $V_1$ and $V_2$.
The edge assignments sum to $q$ at each node in $V_1$ and to 
$p$ at each node in $V_2$ (see, e.g., \cite{bolkertransportation}).
If $p$ and $q$ are distinct primes, the vertices of $P(p,q)$ 
correspond to spanning trees of $K_{p,q}$ 
that satisfy the same conditions on the edge assignments.
The edges of such a spanning tree
are in bijection with the positive entries of the vertex ${\bf x}$.
Hence there are exactly $p+q-1$ such positive entries and 
exactly $(p-1)(q-1)$ zero entries. Using the
bijection established in the above proof we conclude that
each facet of $\C_{pq}$ has exactly $(p-1)(q-1)$ vertices.
Since the dimension of $\C_{pq}$ is equal to 
$(p-1)(q-1)$, each facet must be a simplex.
\end{proof}

%%%%%%%%%%%%%%%%%%%%%%%%%%%%%%%%%%%%%%%%%%%%%%%%%%%%%%%%%%%%%%%%%%%%%%%%%%%

\section{Explicit Computations and a Conjectural Counterexample}\label{expcompsection}
% Used to be called: The cases $m=p$ and $m=2p$ for a prime $p$

We start with the case of $m=p$, a prime. 
The vertices of the cyclotomic polytope $\C_p \subset \R^{ p-1 }$ are the unit vectors $e_1, e_2, \dots, e_{p-1}$, and $-\1 = - \sum_j e_j$. (We apologize for the change of notation from Section \ref{cycpolsection}.)
This simplex has the apparent unimodular triangulation
\[
  \left\{ \conv \left( \0, e_1, e_2, \dots, e_{p-1} \right) , \conv \left( \0, e_1, \dots, e_{ p-2 } , -\1 \right) , \dots, \conv \left( \0, e_2, \dots, e_{p-1} , -\1 \right) \right\} ,
\]
the only triangulation that uses the origin.
Hence Kl\o ve--Parker's Theorem \ref{coordinatorprime} is an immediate consequence of Corollary \ref{maincor} and the fact that the $h$-vector of a simplex is $(1, 1, \dots, 1)$:
\[
  G (x) = \frac{ h_{ \C_p } (x) }{ (1-x)^p }  = \frac{ x^{ p-1 } + x^{ p-2 } + \dots + 1 }{ (1-x)^p } \ .
\]

The second case is $m=2p$ for an odd prime $p$. 
$\C_{ 2p }$ is totally unimodular, and so the facets are supported by hyperplanes of the form
\begin{equation}\label{facethypequation}
  a_1 x_1 + \dots + a_{ p-1 } x_{ p-1 } = 1 \ .
\end{equation}
Furthermore, since
$
  A_{ 2p } = \left[ I_{ p-1 } \ -\!\1 \ -\!I_{ p-1 } \ \1 \right]
$
the $a_j$'s are all $0$ or $\pm 1$.
Let us call two vertices $v_1, v_2$ of $\C_{ 2p }$ \emph{opposite} if $v_2 = -v_1$.
A facet cannot contain two opposite vertices because otherwise the right-hand side of \eqref{facethypequation} would be $0$.

\begin{proposition}\label{2p-ksubsetface}
Suppose $p$ is an odd prime and $k \le \frac{ p-1 }{ 2 }$.
Then every $k$-subset of $A_{ 2p } = \left[ I_{ p-1 } \ -\!\1 \ -\!I_{ p-1 } \ \1 \right]$ that does not contain opposite vertices forms a $(k-1)$-face of $\C_{ 2p }$.
\end{proposition}

\begin{corollary}\label{2p-ksubsetfacenumber}
Suppose $p$ is an odd prime and $k \le \frac{ p-1 }{ 2 }$.
Then $\C_{ 2p }$ has
\[
  f_{ k-1 } = 2^k \binom p k
\]
$(k-1)$-faces.
\end{corollary}

\begin{proof}[Proof of Proposition \ref{2p-ksubsetface}]
Given a $k$-subset $S \subseteq A_{ 2p }$ without opposite vectors, we consider two cases, depending whether or not $\pm \1 \in S$.

\vspace{10pt} \noindent
1.\ case: $\pm \1 \notin S$.

First suppose $k \le \frac{ p-1 }{ 2 } - 1$.
We choose $n$ vectors from $I_{ p-1 }$ and $m = k-n$ vectors from $-I_{ p-1 }$. Without loss of generality, suppose these vectors are $e_1, \dots, e_n, -e_{ n+1 }, \dots, -e_k$.
Set $b = \frac{ m-n }{ p-1-k }$; note that $|b| < 1$ because $k \le \frac{ p-1 }{ 2 } - 1$.
Consider the hyperplane
\[
  x_1 + \dots + x_n - x_{ n+1 } - \dots - x_k + b \left( x_{ k+1 } + \dots + x_{ p-1 } \right) = 1 \ .
\]
Our $k$ chosen vectors are on this hyperplane, and we claim that the remaining vectors in $A_{ 2p }$ satisfy
\begin{equation}\label{offhyperplane}
  x_1 + \dots + x_n - x_{ n+1 } - \dots - x_k + b \left( x_{ k+1 } + \dots + x_{ p-1 } \right) < 1 \ .
\end{equation}
For the remaining unit vectors this follows from $|b| < 1$, and for $x = \pm \1$ \eqref{offhyperplane} becomes the inequality $0<1$.

Now suppose $k = \frac{ p-1 }{ 2 }$.
Again we choose $n$ vectors from $I_{ p-1 }$ and $m = k-n$ vectors from $-I_{ p-1 }$.
If $n \ne 0$ or $k$, the above proof goes through verbatim.
If $n = k$, set $b = -1 + \frac{ 1 }{ 2k }$ and continue the proof above.
If $n = 0$, set $b = 1 - \frac{ 1 }{ 2k }$ and continue the proof above.

\vspace{10pt} \noindent
2.\ case: $\1 \in S$. (The case $-\1 \in S$ is analogous, so that we will omit it here.)

Again we choose $n$ vectors from $I_{ p-1 }$ and $m = k-1-n$ vectors from $-I_{ p-1 }$.
We may assume these vectors are $e_1, \dots, e_n, -e_{ n+1 }, \dots, -e_{k-1}$.
Set $b = \frac{ m-n+1 }{ p-k } $; note that $|b| \le \frac{ k }{ p-k } < 1$.
Consider the hyperplane
\[
  x_1 + \dots + x_n - x_{ n+1 } - \dots - x_{k-1} + b \left( x_k + \dots + x_{ p-1 } \right) = 1 \ .
\]
Our $k$ chosen vectors are on this hyperplane, and again one can easily check that the remaining vectors in $A_{ 2p }$ satisfy
\[
  x_1 + \dots + x_n - x_{ n+1 } - \dots - x_{k-1} + b \left( x_k + \dots + x_{ p-1 } \right) < 1 \ .
\]
\end{proof}

\begin{remark} \rm One can use the correspondence of the facets of $\C_{ 2p }$ to 
the vertices of $P(2,p)$ described in the proof of Proposition \ref{chapman}
to show that $\C_{ 2p }$ has $p \binom{ p-1 }{ \frac{ p-1 }{ 2 } }$ facets.
We do not know the number of facets for the more general cyclotomic polytopes 
$\C_{ pq }$ for distinct primes $p$ and $q$; it would be interesting if the 
correspondence to transportation polytopes could lead to this number.
\end{remark}

Proposition \ref{chapman} and Corollary \ref{2p-ksubsetfacenumber} allow us to prove Parker's Conjecture \ref{parkerconj2}:

\begin{theorem}
The coordinator polynomial of $\Z [\zeta_{2p}]$, where $p$ is an 
odd prime, equals
\[
  h_{ 2p } (x) = \sum_{ k=0 }^{ \frac{ p-3 }{ 2 } } \left( x^k + x^{ p-1-k } \right) \sum_{ j=0 }^{ k } \binom p j + x^{ \frac{ p-1 }{ 2 } } \sum_{ j=0 }^{ \frac{ p-1 }{ 2 } } \binom p j .
\]
\end{theorem}

\begin{proof}
The cyclotomic polytope $\C_{ 2p }$ is simplicial by Proposition \ref{chapman}, so Corollary \ref{maincor} applies.
For $j \le \frac{ p-1 }{ 2 }$, Corollary \ref{2p-ksubsetfacenumber} gives
\[
  h_j 
  = \sum_{ k=0 }^j (-1)^{ j-k } \binom{ p-1-k }{ j-k } f_{ k-1 }
  = \sum_{ k=0 }^j (-1)^{ j-k } \binom{ p-1-k }{ j-k } 2^k \binom p k
  = \sum_{ k=0 }^j \binom p k ,
\] 
as one easily checks that
\[
  \sum_{ k=0 }^j (-1)^{ j-k } \binom{ p-1-k }{ j-k } 2^k \binom p k
  - \sum_{ k=0 }^{j-1} (-1)^{ j-1-k } \binom{ p-1-k }{ j-1-k } 2^k \binom p k
  = \binom p j .
\]
Palindromy of the $h$-vector gives $h_j$ for $j > \frac{ p-1 }{ 2 }$.
\end{proof}

Going beyond $m = p$ or $2p$, next we prove Parker's Conjecture \ref{parkerconj3}.

\begin{corollary}
The coordinator polynomial of $\Z \left[ \zeta_{15} \right]$ equals
\[
  c_{ \Z \left[ \zeta_{15} \right] } (x) = \left( 1 + x^8 \right) + 7 \left( x + x^7 \right) + 28 \left( x^2 + x^6 \right) + 79 \left( x^3 + x^5 \right) + 130 x^4 .
\]
\end{corollary}

\begin{proof}
By Proposition \ref{cyc-product}, the polytope $\C_{ 15 }$ has vertices
\[
  A_{ 15 } =
  \left[ \begin{array}{ccccccccccccccc}
  I_4 & -\1 &     &     & -I_4 & \1 \\
      &     & I_4 & -\1 & -I_4 & \1
  \end{array} \right] ,
\]
and it is simplicial by Proposition \ref{chapman}.
With this data, one can easily use the software {\tt polymake} \cite{polymake} to check that $\C_{ 15 }$ has the $h$-polynomial $x^8 + 7 x^7 + 28 x^6 + 79 x^5 + 130 x^4 + 79 x^3 + 28 x^2 + 7 x + 1$. The result now follows with Corollary \ref{maincor}.
\end{proof}

\begin{figure}[h] 
\begin{tabular}{r|l}
$m$ & $h_{ \Z \left[ \zeta_m \right] }$ \\
\hline
% 2 & $x+1$ \\
% 3 & $x^2 + x + 1$ \\
% 4 & $(x+1)^2$ \\
% 5 & $ x^4 + x^3 + x^2 + x + 1$ \\
 6 & $ x^2 + 4x + 1 $ \\
% 7 & $ x^6 + x^5 + x^4 + x^3 + x^2 + x + 1$ \\
% 8 & $ (x+1)^4 $ \\
% 9 & $ \left( x^2 + x + 1 \right)^3 $ \\
10 & $ x^4 + 6 x^3 + 16 x^2 + 6 x + 1 $ \\
%11 & $ x^{ 10 } + x^9 + x^8 + x^7 + x^6 + x^5 + x^4 + x^3 + x^2 + x + 1$ \\
12 & $ \left( x^2 + 4x + 1 \right)^2 $ \\
%13 & $ x^{ 12 } + x^{ 11 } + x^{ 10 } + x^9 + x^8 + x^7 + x^6 + x^5 + x^4 + x^3 + x^2 + x + 1$ \\
14 & $ x^6 + 8 x^5 + 29 x^4 + 64 x^3 + 29 x^2 + 8 x + 1 $ \\
15 & $x^8 + 7 x^7 + 28 x^6 + 79 x^5 + 130 x^4 + 79 x^3 + 28 x^2 + 7 x + 1$ \\
%16 & $ (x+1)^8 $ \\
%17 & $ x^{ 16 } + x^{ 15 } + \dots + 1$ \\
18 & $ \left( x^2 + 4x + 1 \right)^3 $ \\
%19 & $ x^{ 18 } + x^{ 17 } + \dots + 1$ \\
20 & $ \left( x^4 + 6 x^3 + 16 x^2 + 6 x + 1 \right)^2 $ \\
21 & $ x^{12} + 9x^{11} + 45x^{10} + 158x^9 + 432x^8 + 909x^7 + 1302x^6 + 
\dots + 1 $ \\
22 & $ x^{10} + 12x^9 + 67x^8 + 232x^7 + 562x^6 + 1024 x^5 + 562x^4 + 232x^3 
+ 67 x^2 + 12 x + 1 $ \\
%23 & $ x^{22} + x^{21} + \dots + 1 $ \\
24 & $ \left( x^2 + 4x + 1 \right)^4 $ \\
%25 & $ \left( x^4 + x^3 + x^2 + x + 1 \right)^5 $ \\
26 & $ x^{12} + 14 x^{11} + 92x^{10} + 378x^9 + 1093x^8 + 2380x^7 + 4096x^6 
+ \dots + 1$ \\
%27 & $ \left( x^2 + x + 1 \right)^9 $ \\
28 & $ \left( x^6 + 8 x^5 + 29 x^4 + 64 x^3 + 29 x^2 + 8 x + 1\right)^2 $ \\
%29 & $ x^{28} + x^{27} + \dots + 1 $ \\ 
30 & $  x^8 + 22x^7 + 208x^6 + 874x^5 + 1480x^4 + 874x^3 + 208x^2 + 22x + 1 $ \\
%31 & $ x^{30} + x^{29} + \dots  + 1 $ \\
%32 & $ \left( x+1\right)^{16} $ \\
33 & $ x^{20} + 13x^{19} + 91x^{18} + 444x^{17} + 1677x^{16} + 5187x^{15} + 
13614x^{14} + 31083x^{13} + 61422x^{12} $ \\
   & $ \qquad + 100561x^{11} + 126214x^{10} + \dots + 1 $\\
34 & $x^{16} + 18 x^{15} + 154 x^{14} + 834 x^{13} + 3214 x^{12} + 9402 
x^{11} + 21778 x^{10} + 41226x^9 + 65536 x^8 + \dots + 1 $ \\
35 & $x^{24} + 11x^{23} + 66x^{22} + 286x^{21} + 1001x^{20} + 2996x^{19} + 
7896x^{18} + 18631x^{17} + 39671x^{16} $ \\
   & $ \qquad + 76046x^{15} + 128726x^{14} + 185206x^{13} + 212926x^{12} + \dots + 1 $ \\
36 & $ \left( x^2 + 4x + 1 \right)^6 $ \\
%37 & $ x^{36} + x^{35} + \dots + 1 $ \\
38 & $ x^{18} + 20x^{17} + 191x^{16} + 1160x^{15} + 5036x^{14} + 16664x^{13} 
+ 43796 x^{12} + 94184x^{11} + 169766 x^{10} $ \\
   & $ \qquad + 262144 x^9 + \dots + 1 $  \\
39 & $x^{24} + 15x^{23} + 120x^{22} + 667x^{21} + 2865x^{20} + 10068x^{19} + 
29998x^{18} + 77670x^{17} $ \\ 
   & $ \qquad + 177966x^{16} + 363919x^{15} + 655692x^{14} + 1001649x^{13} + 
1214590x^{12} + \dots + 1 $\\
40 & $ \left( x^4 + 6 x^3 + 16 x^2 + 6 x + 1\right)^4 $ \\
%41 & $x^{40} + x^{39} + \dots + 1 $
\end{tabular}
\caption{The coordinator polynomials of $\Z\left[ \zeta_m\right]$ for $m \leq 41$.} \label{coordinator-table}
\end{figure}

For reference, we give here the first $41$ coordinator polynomials. 
(In the table, we omitted $h_m(x)$ for prime powers $m$.)
The results of this article could be used to compute the coordinator
polynomials for $m \leq 104$, where $m=105$ is the first
non-trivial case that our results do not cover. Although in principal
these coordinator polynomials could be computed, the feasible
range seems to end with $m=41$ with the current 
computational tools like {\tt 4ti2} \cite{4ti2} which we used for 
toric initial ideal computations, and {\tt CoCoA} \cite{cocoa} 
which we used for Hilbert series computations. 
We offer the following conjecture for the first non-trivial case.

\begin{conjecture}
The coordinator polynomial $h_{ 105 }(x)$ is not palindromic.
\end{conjecture}

We conclude by giving some supporting evidence for this conjecture.
We call a polytope $\P \subset \R^d$ \emph{integral} is all its vertices
are in $\Z^d$.
An integral polytope $\P := \left\{ x \in \R^d : \ A x \le \1 \right\}$
that contains the origin in its interior is called \emph{reflexive} if
$A$ is an integral matrix \cite{batyrevdualpoly}.
The \emph{Ehrhart series} of an integral $d$-polytope $\P$ is the rational
generating function $\sum_{ k \ge 0 } \# \left( k\P \cap \Z^d \right) x^k$,
which is of the form $\frac{ f(x) }{ (1-x)^{ d+1 } }$ for some polynomial 
$f$ of degree at most $d$ \cite{ehrhartpolynomial}.
Hibi \cite{hibidual} proved that an integral polytope that contains the 
origin in its interior is reflexive if and only if the numerator of its
Ehrhart series is palindromic.
Since we proved that $\C_m$ is normal if $m$ is divisible by at most two
odd primes, Hibi's theorem implies that $\C_m$ is reflexive for these $m$.

On the other hand, Seth Sullivant \cite{seth} computed some of the facets
of $\C_{ 105 }$ and found that the defining matrix $A$ is not integral,
that is, $\C_{ 105 }$ is \emph{not} reflexive. If one could show that
$\C_{ 105 }$ is normal, then this would provide a counterexample to 
Parker's Conjecture \ref{parkerconj1}.

{\bf Acknowledgements}
We thank Andrew Beyer for help with {\tt polymake} computations, Robin Chapman for the proof of Proposition \ref{chapman}, Vic Reiner for bringing our attention to \cite{martinreinercyclotomic}, G\"unter Ziegler for pointing us to references for the direct sum of polytopes, and Seth Sullivant and an anonymous referee for many helpful comments on an earlier version of this paper.

%%%%%%%%%%%%%%%%%%%%%%%%%%%%%%%%%%%%%%%%%%%%%%%%%%%%%%%%%%%%%%%%%%%%%%%%%%%

%\vfill
%\newpage

%\bibliographystyle{amsplain}
%\bibliography{cyclotomic}

\providecommand{\bysame}{\leavevmode\hbox to3em{\hrulefill}\thinspace}
\providecommand{\MR}{\relax\ifhmode\unskip\space\fi MR }
% \MRhref is called by the amsart/book/proc definition of \MR.
\providecommand{\MRhref}[2]{%
  \href{http://www.ams.org/mathscinet-getitem?mr=#1}{#2}
}
\providecommand{\href}[2]{#2}

\setlength{\parskip}{0cm} 
\end{document}